\documentclass[11pt]{smfart}
\usepackage{color}
\usepackage{amssymb,verbatim}
\usepackage{hyperref}
\usepackage{amsthm,array,amssymb,amscd,amsfonts,amssymb,latexsym, url}
\usepackage{amsmath}
\usepackage[all]{xy}
\usepackage[french]{babel}

\usepackage{cyrillic}

\newcounter{spec}
{\end{list}}

\renewcommand{\P}{{\mathbf P}}
\usepackage{cyrillic}

\newcommand{\Q}{{\mathbb Q}}

\renewcommand{\L}{{\mathbb L}}

\renewcommand{\lim}{\varprojlim}

\numberwithin{equation}{section}

\newfont{\gothic}{eufb10}

\renewcommand{\qed}{{\hfill$\square$}}

\newtheorem{theo}{Th\'{e}or\`{e}me}[section]
\newtheorem{prop}[theo]{Proposition}

\newtheorem{lem}[theo]{Lemme}
\newtheorem{cor}[theo]{Corollaire}
\theoremstyle{definition}
\newtheorem{defi}[theo]{D\'efinition}
\theoremstyle{remark}
\newtheorem{rema}[theo]{Remarque}

\newcommand{\bthe}{\begin{theo}}
\newcommand{\ble}{\begin{lem}}
\newcommand{\bpr}{\begin{prop}}
\newcommand{\bco}{\begin{cor}}
\newcommand{\bde}{\begin{defi}}
\newcommand{\ethe}{\end{theo}}
\newcommand{\ele}{\end{lem}}
\newcommand{\epr}{\end{prop}}
\newcommand{\eco}{\end{cor}}
\newcommand{\ede}{\end{defi}}

\DeclareFontFamily{U}{wncy}{}
\DeclareFontShape{U}{wncy}{m}{n}{%
<5>wncyr5%
<6>wncyr6%
<7>wncyr7%
<8>wncyr8%
<9>wncyr9%
<10>wncyr10%
<11>wncyr10%
<12>wncyr6%
<14>wncyr7%
<17>wncyr8%
<20>wncyr10%
<25>wncyr10}{}
\DeclareMathAlphabet{\cyr}{U}{wncy}{m}{n}

\begin{document}

 \title[Espaces lin\'eaires sur les intersections de deux quadriques]
 {Principe local global pour les espaces lin\'eaires sur les intersections de deux quadriques}

\author{J.-L. Colliot-Th\'el\`ene}
\address{C.N.R.S., Universit\'e Paris Sud\\Math\'ematiques, B\^atiment 425\\91405 Orsay Cedex\\France}
\email{jlct@math.u-psud.fr}
  
\date{19 novembre 2014}
\maketitle

 \begin{abstract}
 Dans \url{http://arxiv.org/abs/1410.5671},  Jahnel et Loughran \'etablissent un principe local-global pour
l'existence d'espaces lin\'eaires de dimension $r$ dans les intersections
compl\`etes 
lisses de deux quadriques dans un espace projectif de dimension $2r+2$.
On donne une d\'emonstration alternative d'un r\'esultat un peu plus g\'en\'eral.
  \end{abstract}
 
 \begin{altabstract} In \url{http://arxiv.org/abs/1410.5671}, Jahnel and Loughran
 prove the local principle for linear spaces of dimension $r$ on smooth complete
 intersections of two quadrics in projective space of dimension $2r+2$.
 We present an alternative proof of a slightly more general result.
 \end{altabstract}
 
 \section{Le th\'eor\`eme}
 
\begin{theo}\label{principal}
Soit $k$ un corps de nombres. Soient $n$ un entier naturel et $f(x_{0}, \dots, x_{n})$
et $g(x_{0},\dots, x_{n})$ deux formes quadratiques \`a coefficients dans $k$.
Supposons qu'il existe une forme non singuli\`ere dans le pinceau de formes quadratiques
$\lambda f+\mu g$.  Dans chacun des cas suivants :

(a) $n=2r+2$,

(b) $n=2r+1$,

\noindent la $k$-vari\'et\'e $X$ d\'efinie par $f=g=0$ dans $\P^n_{k}$ contient un sous-espace lin\'eaire
d\'efini sur $k$ de dimension $r$ si et seulement si elle en poss\`ede un sur chaque compl\'et\'e $k_{v}$ de $k$, pour
$v$ parcourant l'ensemble $\Omega$ des places de $k$.

\end{theo}

\begin{proof}

Soit $t$ une variable. La forme quadratique $q(t):= f+tg$ sur le corps $k(t)$ est par hypoth\`ese non d\'eg\'en\'er\'ee. Soit $\delta (t) \in k(t)^{\times}$ 
son discriminant.

D'apr\`es le th\'eor\`eme d'Amer \cite{amer, Pf, L}, sur tout corps $F$ contenant $k$,
la $F$-vari\'et\'e  $X\times_{k}F$ contient un $F$-sous-espace lin\'eaire de dimension $r$
si et seulement si la quadrique dans $\P^n_{F(t)}$ d\'efinie par la
 forme quadratique $q(t)= f+tg$ sur le corps $F(t)$ s'annule 
 sur un espace lin\'eaire de dimension $r$.
 Comme la forme quadratique  $q(t)$ est   non d\'eg\'en\'er\'ee, cette
derni\`ere propri\'et\'e \'equivaut
 au fait que
 la forme quadratique $q(t)$ 
 contient en facteur direct la somme orthogonale de $r+1$ facteurs hyperboliques $H=<1,-1>$.
 
 Pour $n=2r+2$, ceci \'equivaut au fait que la forme quadratique 
 $f+tg$ est
 isomorphe \`a la forme quadratique $q_{1}(t):= (r+1).H \perp <(-1)^{r+1}.\delta(t)>.$
 
 Pour $n=2r+1$, ceci est \'equivalent au fait  que la  forme quadratique 
  $f+tg$
 est isomorphe \`a la forme quadratique $q_{2}(t):= (r+1).H$.

 D'apr\`es \cite[Prop. 1.1]{CTCS},  l'application diagonale de groupes de Witt
 $$W(k(t)) \to \prod_{v \in \Omega} W(k_{v}(t))$$
 est injective.
 
 Le th\'eor\`eme de simplification de Witt permet de d\'eduire de ce r\'esultat l'\'enonc\'e :
 Si deux  formes quadratiques sur $k(t)$ sont isomorphes sur chacun des $k_{v}(t)$ pour $v \in \Omega$,
 alors elles sont isomorphes sur $k(t)$.
 
 Ceci ach\`eve la d\'emonstration.
 \end{proof}
 
 \bigskip
 
 \section{Remarques}
 
\begin{rema}
 
 Dans [JL],  Thm. 1.5, Jahnel et Loughran consid\`erent le cas o\`u $X$ est une intersection compl\`ete lisse
 de deux quadriques et o\`u  $n=2r+2$. Ils montrent dans ce cas le r\'esultat plus pr\'ecis : si pour presque toute place $v$ de $k$,
 la $k$-vari\'et\'e $X\times_{k}k_{v}$ contient un $k_{v}$-espace lin\'eaire de dimension~$r$,
 alors elle contient un $k$-espace lin\'eaire de dimension $r$. Ceci tient au fait que dans leur cas
 le principe local-global qui est utilis\'e est : un \'el\'ement d'un corps de nombres $K$ est un carr\'e
 si et seulement si il l'est sur presque tout compl\'et\'e de $K$.  
 
  \end{rema}
  
  \begin{rema}
  Dans le cas $n=2r+1$, si $X$ est une intersection compl\`ete lisse de deux quadriques,
  la dimension maximale d'un sous-espace lin\'eaire contenu dans $X$ est $r-1$
  (\cite[Cor. 2.4]{reid}). Le th\'eor\`eme est donc vide dans ce cas. Il n'est int\'eressant que
  pour $X$ singuli\`ere.

    \end{rema}
    
      \begin{rema} Comme me le fait observer David Leep (20 f\'evrier 2015),  dans le th\'eor\`eme~\ref{principal} on ne peut pas omettre l'hypoth\`ese de non singularit\'e
      g\'en\'erique dans le pinceau. Soit $C \subset \P^3_{\Q}$ une courbe de genre 1 contre-exemple au principe de Hasse, donn\'ee
      par un syst\`eme de deux formes quadratiques $$q_1(x_1, x_2, x_3, x_4)= q_2(x_1, x_2, x_3, x_4)=0.$$

Les deux formes en 7 variables $Q_{1}=q_1(x_1, \ldots, x_4) + x_5 x_6$ et $Q_{2}=q_2(x_1, \ldots, x_4) + x_5 x_7$
      s'annulent simultan\'ement sur un espace vectoriel de dimension 3 sur chaque compl\'et\'e de $\Q$ mais
      pas sur $\Q$.  On le voit en appliquant  le th\'eor\`eme d'Amer-Brumer \cite{amer, L}
      \`a $Q_{1}+tQ_{2}$. Ceci correspond au cas  $r=2$ et  $n=2r + 2 = 6$ du th\'eor\`eme  \ref{principal}.
      
      De m\^eme, les deux formes en 10 variables $Q_{1}=q_1(x_1, \ldots, x_4) + x_5 x_6 + x_8 x_9$ et $Q_{2}=q_2(x_1, \ldots, x_4) + x_5x_7 + x_8 x_{10}$ s'annulent sur un espace vectoriel de dimension 5 sur chaque compl\'et\'e de $\Q$, mais pas sur $\Q$.
      Ceci correspond au cas $ r=4$ et $ n=2r + 1 = 9$ du th\'eor\`eme  \ref{principal}.
        \end{rema}


\begin{thebibliography}{99}


\bibitem{amer} M. Amer, Quadratische Formen \"uber Funktionenk\"orpern, Dissertation, Johannes Gutenberg Universit\"at, 
Mainz 1976.
 
 
\bibitem{CTCS}   J.-L. Colliot-Th\'el\`ene, D. Coray et J.-J. Sansuc, Descente et principe de Hasse pour certaines vari\'et\'es rationnelles,
J. reine und angew. Math. {\bf  320} (1980) 150--191.
 




\bibitem{JL} J. Jahnel et D. Loughran, The Hasse principle for lines on Del Pezzo surfaces, 
\url{http://arxiv.org/abs/1410.5671}

\bibitem{L} D. Leep, The Amer--Brumer theorem over arbitrary fields
 \url{http://www.ms.uky.edu/~leep/Amer-Brumer_theorem.pdf}

\bibitem{Pf} A. Pfister, Quadratic Forms with Applications to Algebraic Geometry and Topology,
Cambridge University Press (1995).

\bibitem{reid} M. Reid, The complete intersection of two or more quadrics, Ph.D. thesis, Cambridge, 1972.

\end{thebibliography}
\end{document}